\documentclass[12pt]{article}
\usepackage{amssymb,amsmath}
\numberwithin{equation}{section}

\def\DD{\hbox{$I\hskip -4pt D$}}

\def\Z{\mathbb{Z}}

\def\C{\mathbb{C}}

\def\qed{\hfill\ \ifhmode\unskip\nobreak\fi\ifmmode\ifinner
 \else\hskip5pt\fi\fi \hbox{\hskip5pt\vrule width4pt height6pt depth1.5pt
 \hskip 1 pt}}
\def\a{\alpha}

\def\dsum#1#2{\mbox{$\sum\limits_{#1}^{#2}$}}
\def\b{\beta}
\def\d{\delta}

\def\g{\gamma}
\def\G{\Gamma}
\def\l{\lambda}

\def\dt{\mbox{$\frac{d}{dt}$}}
\def\Vir{\mbox{Vir}}

\def\si{\sigma}
\def\sc{\scriptstyle}
\def\ssc{\scriptscriptstyle}
\def\F{\mathbb{F}}

\def\cl{\centerline}
\def\DD{{\cal D}}

\def\ol{\overline}

\def\wt{\widetilde}
\def\wh{\widehat}
\def\rar{\rightarrow}

\def\bs{\backslash}
\def\hs{\hspace*}

\def\vs{\vspace*}
\def\ra{\rangle}
\def\la{\langle}
\def\ni{\noindent}
\def\hi{\hangindent}
\def\ha{\hangafter}
\def\WW{{\cal W}}
\def\AA{{\cal A}}

\def\nb{}%#1{\mbox{#1}}
\begin{document}

\cl {{\Large Classification of quasifinite
$\WW_\infty$-modules}\footnote{Supported by NSF grants 10471096,
10571120 of China and ``One Hundred Talents Program'' from
University of Science and Technology of China}} \vs{4pt}

\cl{(to appear in {\it Israel Journal of Mathematics})}\vs{10pt}
 \cl{Yucai Su$^*$ \ \ \ and \ \ \ Bin Xin$^\dag$} \vs{4pt}

  \cl{\small
$^{*}$Department of Mathematics, University of Science and
Technology of \vs{-2pt}China} \cl{\small Hefei 230026, China}
\vs{6pt}

\cl{\small $^{\dag}$Department of Mathematics, Shanghai Jiaotong
University, %} \cl{\small
Shanghai 200240, China}\vs{6pt}

 \cl{\small E-mail: ycsu@sjtu.edu.cn,
xinbinsu@sjtu.edu.cn}
\vs{15pt}

\ni{\small {\bf Abstract.}
It is proved that an irreducible quasifinite $\WW_\infty$-module is a
highest or lowest weight module or a module of the
intermediate series; a uniformly bounded
indecomposable weight \nb{$\WW_\infty$-module} is a
module of the intermediate
series.
For a nondegenerate additive \nb{subgroup} $\G$ of $\F^n$, where $\F$ is a
field of characteristic zero, there is a \nb{simple}
Lie or \nb{associative} algebra $\WW(\G,n)^{(1)}$
spanned by differential operators
$uD_1^{m_1}\cdots D_n^{m_n}$ for $u\in\F[\G]$ (the
group algebra), and $m_i\ge0$ with
$\sum_{i=1}^n m_i\ge1$, where $D_i$ are degree
operators.
It is also proved that an indecomposable
quasifinite weight $\WW(\G,n)^{(1)}$-module is a module
of the intermediate series if $\G$ is not isomorphic to
$\Z$.
\vs{5pt}

\ni
{\it
Mathematics Subject Classification
(1991):} 17B10, 17B65, 17B66, 17B68
\vs{15pt}}

\ni{\bf1. Introduction.}
Let us start with the general definition.
For an algebraically closed field $\F$ of
characteristic zero, let $\G$ be a {\bf nondegenerate}
additive subgroup of $\F^n$, i.e., it contains an
$\F$-basis of $\F^n$.
Let $\F[\G]={\rm
span}\{t^\a\,|\, \a \in \G\}$ denote the group
algebra of $\G$ with the algebraic operation
$t^\a\cdot t^\b = t^{\a+\b}$
for
$\a,\b \in \G$.
We define the {\bf degree operators}
$D_i$ to be the derivations of $\F[\G]$ determined by
$D_i:t^\a \mapsto \a_it^\a$
for
$\a \in \G,\,i = 1,...,n$.
Here and below, an element
$\a\in\F^n$ is always written as
$\a = (\a_1,...,\a_n)$.
The {\bf Lie algebra
$\WW(\G,n)$ of Weyl type} [S4] is a tensor product space
of the group algebra $\F[\G]$ with the polynomial
algebra $\F[D_1,...,D_n]$:
$$
\WW(\G,n)=\F[\G]\otimes\F[D_1,...,D_n]={\rm span}
\{t^\a D^\mu\,|\,\a\in\G,\,\mu\in\Z_+^n\},
\eqno(1.1)
$$
where $D^\mu=\prod_{i=1}^{n}D_i^{\mu_i},$
with the Lie bracket:
$$
[t^\a D^\mu,t^\b D^\nu]= (t^\a D^\mu)\cdot(t^\b
D^\nu) -(t^\b D^\nu) \cdot (t^\a D^\mu),
$$
and
$$
(t^\a D^\mu)\cdot(t^\b
D^\nu)=\sum_{\l\in\Z_+^n}
\left(\begin{array}{c}\mu\\\l\end{array}\right)
\b^\l t^{\a+\b}D^{\mu+\nu-\l},
\eqno(1.2)
$$
where
$\b^\l=\prod_{i=1}^n\b_i^{\l_i}$ (here without
confusion, we use notation $\b^\l$ similar to notation
$D^\mu$ in (1.1)), and
$(^\mu_\l) = \prod_{i=1}^n(^{\mu_i}_{\l_i})$. Furthermore, for
$i,j\in\F$,
$(^i_j) = i(i - 1)\cdots(i - j + 1)/j!$ if
$j\in\Z_+$, or $(^i_j) = 0$ otherwise.

It is proved [S3] that $\WW(\G,n)$ has a nontrivial
universal central extension if
and only if $n=1$.
The Lie bracket for the universal central extension
$\wh\WW(\G,1)$ of
$\WW(\G,1)$ is defined
by
$$
\begin{array}{ll}
[t^\a [D]_\mu,t^\b [D]_\nu]
=\!\!\!&
(t^\a [D]_\mu) \cdot (t^\b [D]_\nu)
 - (t^\b [D]_\nu) \cdot (t^\a [D]_\mu)
\\[4pt]
&
 + \ \d_{\a,-\b}(-1)^\mu\mu!\nu!
\left(\begin{array}{c}\a + \mu\\ \mu + \nu + 1\end{array}\right)C,
\end{array}
\eqno(1.3)
$$
for $\a,\b \in \G \subset \F,\,\mu,\nu \in \Z_+,$
where $[D]_\mu = D(D - 1)
\cdots(D - \mu + 1)$, and $C$
is a central element
of $\wh\WW(\G,1)$.
The 2-cocycle of $\WW(\Z,1)$
corresponding to (1.3) seems to appear first in [KP].

Denote by $\WW(\G,n)^{(1)}$ the Lie subalgebra of
$\WW(\G,n)$ spanned by
$\{t^\a D^\mu\,|\,\a\in\G,\,|\mu|\ge1\}$, where
$|\mu|=\sum_{i=1}^n\mu_i$.
Similarly, we can define
$\wh\WW(\G,1)^{(1)}$.
Then
$\WW_{1+\infty}=\wh\WW(\Z,1)$ and
$\WW_\infty=\wh\WW(\Z,1)^{(1)}$ are the well-known
$\WW$-infinity
algebras, which arise naturally in various physical
theories such as conformal field theory, the theory
of the quantum Hall effect,
etc.~and which receive intensive studies in the
literature
(cf.~[BKLY,\,FKRW,\,KL,\,KR1,\,KR2,\,KWY,\,S4]).

Note that $\WW(\G,n)^{(1)}$ is also an associative
algebra under the product (1.2).
It can be proved
that $\WW(\G,n)^{(1)}$ is simple as a Lie or
associative algebra (cf.~[SZ1]).
We denote it
by $\AA(\G,n)^{(1)}$ when we consider it as an
associative algebra.
Clearly an
$\AA(\G,n)^{(1)}$-module is also a
$\WW(\G,n)^{(1)}$-module, but not necessarily the
converse.
Thus
it suffices to consider $\WW(\G,n)^{(1)}$-modules.
The Lie algebra
$\WW(\G,n)^{(1)}= \oplus_{\a\in\G}\WW(\G,n)^{(1)}_\a$
is $\G$-graded with the grading space
$$
\WW(\G,n)^{(1)}_\a={\rm span}\{t^\a
D^\mu\,|\,\mu\in\Z_+^n\bs\{0\}\}
\mbox{ \ \ for \ \ }\a\in\G.
\eqno(1.4)
$$
In [S4], one of us classified
the quasifinite
modules over $\WW(\G,n)$.
In this paper, we shall
consider the
more difficult problem of classifying the
quasifinite modules over
$\WW(\G,n)^{(1)}$.
Here, a $\WW(\G,n)^{(1)}$-module
$V$ is called
a {\bf quasifinite module} if $V=\oplus_{\a\in\G}V_\a$
is a $\G$-graded $\F$-vector space such that
$\WW(\G,n)^{(1)}_\a V_\b\subset V_{\a+\b},
\,{\rm dim}V_\a<\infty$ for $\a,\b\in\G$.
When we study the representations of Lie algebras
of this kind, since each grading space in (1.4) is
still infinite-dimensional,
the classification of quasifinite modules is thus a
nontrivial problem, as pointed in [KL].

For $\a\in\F^n$, one can define quasifinite
$\WW(\G,n)^{(1)}$- or
$\wh\WW(\G,1)^{(1)}$-modules $A_\a,B_\a$ as follows:
They have
basis $\{y_\b\,|\,\b\in\G\}$ such that the
central element $C$
acts trivially and
$$
\begin{array}{ll}
A_\a:&
(t^\b D^\mu)y_\g=(\a+\g)^\mu y_{\b+\g},\vs{4pt}\\
B_\a:&
(t^\b D^\mu)y_\g=(-1)^{|\mu|+1}(\a+\b+\g)^\mu y_{\b+\g},
\end{array}
$$
for $\b,\g \in \G,\,\mu
 \in \Z_+^n\bs\{0\}$ (where
$(\a+\g)^\mu$ is a notation as $\b^\l$ in (1.2)).
These modules
are defined in [S4, Z].
Obviously, $A_\a$ or $B_\a$
is irreducible
if and only if $\a \notin \G$.
Clearly $A_\a$ is
also an $\AA(\G,n)^{(1)}$-module, but not $B_\a$.
We refer any subquotient module of $A_\a$ or $B_\a$
to as a {\bf module of the intermediate series}
(cf.~[S4]).
Then the main result of the present paper
is the following.
\vs{5pt}

{\bf Theorem 1.1.} \
{\it {\rm(i)}
An irreducible
quasifinite module over
\linebreak[4]
$\WW(\Z,1)^{(1)}$ or over
$\WW_\infty = \wh\WW(\Z,1)^{(1)}$ is a highest or
lowest weight module, or a module of the intermediate
series.

{\rm(ii)} An
irreducible quasifinite $\WW(\G,n)^{(1)}$- or
$\wh\WW(\G,1)^{(1)}$-module is a module of the
intermediate series
if $\G$ is not isomorphic to $\Z$.}
\vs{5pt}

Since the complete description of irreducible quasifinite highest
weight modules was obtained in [KL] and lowest weight modules are dual of
highest weight modules, Theorem 1.1 and results in [KL] in fact give a
complete classification
of irreducible quasifinite modules.
Theorem 1.1 also gives a classification of irreducible
quasifinite modules over the associative algebras
$\AA(\G,n)^{(1)}$.

The analogous results to the above
theorem for affine Lie algebras,
the Virasoro algebra, higher rank Virasoro algebras and
Lie algebras of Weyl type or Block type have been
obtained in [C,\,M,\,S4,\,S5,\,S6] (also, cf.~[S2]).

A quasifinite module $V$ is {\bf uniformly bounded}
if there exists $N\ge0$ such that
${\rm dim\,}V_\b\le N$
for all $\b\in\G$;
it is called a {\bf weight module} if $D_1,...,D_n$
are semi-simple operators on $V$.
\vs{5pt}

{\bf Theorem 1.2.} \
{\it {\rm(i)}
A uniformly bounded
indecomposable weight $\WW(\Z,1)^{(1)}$- or
$\WW_\infty$-module is
a module of the intermediate series.

{\rm(ii)}
A quasifinite
indecomposable weight $\WW(\G,n)^{(1)}$- or
$\wh\WW(\G,n)^{(1)}$-module is a module of the
intermediate series
if $\G$ is not isomorphic to $\Z$.}
\vs{5pt}

Finally we would like to point out that although the
main result
of the present paper is similar to that of [S4], one
can see below
that the proof is more technical than that of [S4] due
to the fact
that the elements $t^\b=t^\b D^0,\,\b\in\G,$ do not
appear in
$\WW(\G,n)^{(1)}$.
\vs{15pt}

\ni{\bf2. Quasifinite $\WW_\infty$-modules.}
First we prove Theorem 1.1(i) and Theorem 1.2(i).
We shall only work
on the
non-central extension case since the proof of
the central
extension case is similar.

Now consider the Lie algebra
$W:=\WW(\Z,1)^{(1)}={\rm
span}\{t^iD^j\,|\,i\in\Z,\,j\in\Z_+\bs\{0\}\}$.
In this case
$D=t\frac{d}{dt}$, and by (1.4),
$W=\oplus_{i\in\Z}W_i$
is $\Z$-graded with
$$
W_i={\rm
span}\{t^iD^j\,|\,j\in\Z_+\bs\{0\}\}=
\{t^iDf(D)\,|\,f(D)\in\F[D]\}
$$
for $i\in\Z.$
By (1.2), we have
\setcounter{equation}{0}\setcounter{section}{2}
\begin{eqnarray}
\hs{-40pt}&&[t^iDf(D),t^jDg(D)]
\nonumber\\
\hs{-40pt}&&
=
t^{i+j}D((D+j)f(D+j)g(D)-(D+i)g(D+i)f(D)),
\end{eqnarray}
for $i,j\in\Z,\,f(D),g(D)\in\F[D].$
Also,
$W$ has a triangular decomposition $W=W_+\oplus
W_0\oplus W_-$,
where in general, for any $\Z$-graded space $M$,
we always use
notations $M_+,M_-,M_0$ and $M_{[p,q)}$ to denote the
subspaces
spanned by elements of degree $k$ with $k>0,k<0,k=0$
and $p\le k<q$ respectively.
Denote $\Vir =
 \oplus_{i\in\Z}\F{\sc\,} t^iD$,
which is the (centerless) Virasoro algebra.
\vs{5pt}

{\bf Lemma 2.1.} \
{\it
Let $S$ be a subspace
of $W_0$ with
finite co-dimension.
Given $i_0>0$, let $M_{i_0,S}$
denote the
subalgebra of $W$ generated by
$t^{i_0}D,t^{i_0+1}D,t^{i_0}D^2$
and $S$.
Then there exists some integer $K > 0$
such that
$W_{[K,\infty)} \subset  M_{i_0,S}$.
}
\vs{5pt}

{\it Proof.}
By the
assumption of $S$, there exists some integer
$m_0 \ge 0$ such
that for all integer $m \ge  m_0$, there exists a
polynomial $Df(D) \in  S$ with ${\rm deg\,}f = m$.
We shall prove by induction on $m$ the following claim.
\vs{4pt}

{\bf Claim 1}.
For any $m\in\Z$ with $1 \le  m
\le  m_0$,
there exists some integer $K_m > mK_{m-1}$ (where we
take $K_0=i_0$) such that $t^kD^m \in  M_{i_0,S}$ for
all integers $k \ge  K_m$.
\vs{4pt}

Suppose $m = 1$.
For any integer $k$ sufficiently
large enough, we can write $k=k_1i_0+k_2(i_0 + 1)$
for some $k_1,k_2\in\Z_+\bs\{0\}$, so $t^kD$ can be
generated by $t^{i_0}D,
t^{i_0+1}D$, i.e., $t^kD\in M_{i_0,S}$.
Thus we can
take some integer $K_1>i_0$ large enough to ensure
that the claim holds for $m=1$.
Suppose $1<m\le m_0$
and inductively assume that the claim
holds for $m-1$.
Take $K_m=mK_{m-1}+i_0$.
Then for
any $k\ge K_m$, by (2.1) we have
$a{\sc\,}t^kD^m\equiv
[t^{k-i_0}D^{m-1},t^{i_0}D^2]
\equiv0\ ({\rm mod\,}M_{i_0,S}),$
where $a=((m-1)i_0-2(k-i_0))<0$,
i.e., $t^kD^m\in M_{i_0,S}$.
Thus
the claim holds for $m$.
\vs{5pt}

Now take $K=K_{m_0}$.
For any integer $k\ge K$, we
can now prove
by induction on $m\ge1$ that $t^kD^m\in M_{i_0,S}$
as follows:
If $m\le m_0$, this immediately follows from Claim 1.
Assume that $m>m_0$.
Let $f(D)$ be a polynomial of degree
$m-1\ge m_0$ such
that $Df(D)\in S$, then by (2.1),
$kmt^kD^m\equiv[t^kD,Df(D)]
\equiv0\ ({\rm mod\,}M_{i_0,S})$.
This
proves that $W_{[K,\infty)}\subset M_{i_0,S}$.
\qed
\vs{5pt}

{\bf Lemma 2.2.} \
{\it
Assume that $V$ is an
irreducible quasifinite $W$-module without a highest
or lowest weight.
For any
$i,j\in\Z,\,i\ne0,-1$, the linear map
$$
t^iD|_{V_j}\oplus
t^{i+1}D|_{V_j}\oplus t^iD^2|_{V_j}:\ V_j\to
V_{i+j}\oplus
V_{i+j+1}\oplus V_{i+j}
$$
is injective.
In particular, ${\rm
dim\,}V_j\le 2({\rm dim\,}V_0)+{\rm dim\,}V_1$
for $j\in\Z$.
}
\vs{5pt}

{\it Proof} (cf.~[S4]).
Being irreducible, $V$
must be a weight
module, i.e., there exists some $\a\in\F$, such that
$$
V_i=\{v\in V\,|\,Dv=(\a+i)v\}.
\eqno(2.2)
$$
Say,
$i>0$ and
$(t^iD)v_0=(t^{i+1}D)v_0=(t^iD^2)v_0=0$ for some
$0\ne v_0\in V_j$.
By shifting the grading index of $V_j$
if necessary, we can
suppose $j=0$.
Let $S$ be the kernel of the linear
map $W_0\to{\rm End\,}(V_0): D^m\mapsto D^m|_{V_0}$ for
$m \ge 1$.
Since
${\rm dim\,}V_0 < \infty$, $S$ is a subspace
of $W_0$ with
finite co-dimension.
Then $M_{i,S}v_0=0$ and by
Lemma 2.1, we have
$W_{[K,\infty)}v_0 = 0$ for some $K > 0$.

For any subspace $M$ of $W$, we use $U(M)$ to denote
the subspace,
which is the span of the standard monomials
with respect to a
basis of $M$, of the universal enveloping algebra
of $W$.
Since
$W=W_{[1,K)}+W_0+W_-+W_{[K,\infty)}$, using
the PBW theorem and
the irreducibility of $V$, we have
\setcounter{equation}{2}
\begin{eqnarray}
V
&\!\!\!=\!\!\!&
U(W)v_0=U(W_{[1,K)})U(W_0+W_-)U(W_{[K,\infty)})v_0
\nonumber\\
&\!\!\!=\!\!\!&
U(W_{[1,K)})U(W_0+W_-)v_0.
\end{eqnarray}
Note that $V_+$ is a
$W_+$-module.
Let $V'_+$ be the $W_+$-submodule
of $V_+$ generated
by $V_{[0,K)}$.
We want to prove that $V_+=V'_+.$

So let $k\ge0$ and let $x\in V_+$ have
degree ${\rm deg\,}x=k$.
If
$0\le k<K$, then by definition, $x\in V'_+$.
Suppose $k\ge K$.
Using (2.3), $x$ is a linear combination of the
form $u_1x_1$ with
$u_1\in W_{[1,K)},x_1\in V$.
Thus the
degree ${\rm deg\ssc\,}u_1$
of $u_1$ satisfies $1\le{\rm deg\ssc\,}u_1<K$,
so $0<{\rm
deg\ssc\,}x_1=k-{\rm deg\ssc\,}u_1<k$.
By
inductive hypothesis,
$x_1\in V'_+$, and thus $x\in V'_+$.
This proves that $V_+=V'_+.$

The fact that $V_+=V'_+$ means that the
$W_+$-module $V_+$ is
generated by the finite dimensional space
$V_{[0,K)}$.
Choose a
basis $B$ of $V_{[0,K)}$. Then for any $x\in B$,
we have
$x=u_xv_0$ for some $u_x\in U(W)$.
Regarding $u_x$
as a polynomial
with respect to a basis of $W$, by induction on
the polynomial
degree and using the formula
$[w,w_1w_2] =[w,w_1]w_2+w_1[w,w_2]$
for $w\in W$, $w_1,w_2\in U(W)$, we see that there
exists a positive integer $k_x>K$ sufficiently large
enough such that
$[W_{[k_x,\infty)},u_x]\subset U(W)W_{[K,\infty)}$.
Then from
$W_{[K,\infty)}v_0=0$, we have $W_{[k_x,\infty)}x
=[W_{[k_x,\infty)},u_x]v_0+u_xW_{[k_x,\infty)}v_0
=0$.
Take
$K'={\rm max}\{k_x\,|\,x\in B\}$, then
$W_{[K',\infty)}V_{[0,K)}=0$ and
$$
W_{[K',\infty)}V_+=W_{[K',\infty)}U(W_+)V_{[0,K)}=
U(W_+)W_{[K',\infty)}V_{[0,K)}=0.
$$
Since there exists some integer $K_1>K'$
sufficiently large enough
to ensure that $W_+\subset
W_{[K',\infty)}+[W_{[-K_1,0)},W_{[K',\infty)}]$,
this means that
we have $W_+V_{[K_1,\infty)}=0$.
Now Suppose $x\in
V_{[K_1+K,\infty)}$.
Then by (2.3), it is a sum of
elements of the
form $u_1x_1$ such that $u_1\in W_{[1,K)}$.
But then $x_1$ has
degree ${\rm deg\,}x_1>{\rm deg\,}x-K\ge K_1$,
so $x_1\in
V_{[K_1,\infty)}$.
Thus from
$W_+V_{[K_1,\infty)}=0$, we have
$u_1x_1=0$, i.e., $x=0$.
This proves that $V$
has no degree $\ge K_1+K$.

Now let $K''$ be the maximal integer such that
$V_{K''}\ne0$.
Since $W_0$ is commutative, there exists a
common eigenvector
$v'_0\in V_{K''}$ for $W_0$.
Then $v'_0$ is a
highest weight
vector of $W$, this contradicts the assumption of
the lemma.
\qed
\vs{5pt}

Theorem 1.1(i) will follow from Theorem 1.2(i) and
Lemma 2.2, so it suffices to prove Theorem 1.2(i).
Thus from now on, we suppose
$V$ is a uniformly bounded indecomposable weight
$W$-module such that (2.2) holds.

Regarding
$V$ as a weight module over the Virasoro algebra
$\Vir$, by [S2], there exists some $N\ge0$ such
that ${\rm dim\,}V_k=N$ for all $k\in\Z$ with
$k+\a\ne0$, where $\a \in \F$ is fixed such
that (2.2) holds, and $V$ has only a
finite composition factors as a $\Vir$-module, and
$t^{-1}D|_{V_k}:V_k\rar V_{k-1}$ is bijective when
$k>>0$.
So, we can find a basis
$Y_k=(y_k^{(1)},...,y_k^{(N)})$ of $V_k$ such
that
$$
(t^{-1}D)Y_k=Y_{k-1}\mbox{ \ for \ }k>>0.
\eqno(2.4)
$$
We shall
assume that $N\ge1$ since the proof is trivial if
$N=0$.
In the
following, we always suppose that $k$ is an integer
such that $k>>0$.
Assume that
$$
(t^iD)Y_k=Y_{k+i}P_{i,k}\mbox{ \,for some \,}
N\times N\mbox{
\,matrices \,}P_{i,k} \mbox{ \,and \,}i\in\Z.
$$
By (2.2), (2.4)
and applying $[t^{-1}D,t^iD]=(i+1)t^{i-1}D$ to $Y_k$
for $i=1,2$, we obtain
$$
P_{-1,K}=1,\;\; P_{0,k}=\ol k,\;\; P_{1,k}= [\ol
k]^2+P_1,\;\; P_{2,k}=[\ol k]^3+3\ol kP_1+P_2,
\eqno(2.5)
$$
for
some $N\times N$ matrices $P_1,P_2$.
Here and below,
for convenience, we always identify a scalar $a\in\F$
with the corresponding $N \times  N$ scalar matrix
$a\cdot{\bf1}_N$ when
the context is clear, where ${\bf1}_N$ is the
$N \times  N$ identity matrix.
We also denote
$\ol k = k + \a$ for
$k \in \Z$, and in general, we use the notation
$[a]^j = a(a + 1)\cdots(a + j - 1)$ for
$a\in\F$, $j\in\Z_+$ (cf.~notation $[D]_j$ in
(1.3)).
By choosing a composition series of $V$
regarding as a $\Vir$-module, we can
suppose $P_1,P_2$ are upper-triangular
matrices.
Applying $[tD,t^2D]=t^3D$ to
$Y_k$, by (2.5), we obtain
$$
P_{3,k}=[\ol k]^4+6[\ol k]^2P_1+4\ol kP_2+P_3,
\eqno(2.6)
$$
where $P_3=-3(2P_1+P_1^2-2P_2)+[P_1,P_2],$
and
$[P_1,P_2] = P_1P_2 - P_2P_1$ is the usual
Lie bracket.
Recall that $D = t\dt$.
From this,
one has
$t^{i+j}(\dt)^j = t^i[D]_j$ for
$i \in \Z,\,j \in \Z_+\bs\{0\}$.
In
the following, we shall often use notation
$\dt$ instead of $D$ whenever it is convenient.
Remind that $\dt$ is an operator of degree $-1$.
Assume that
$$
\bigl(\dt\bigr)^iY_k=Y_{k-i}Q_{i,k}\mbox{ \,for some \,}
N\times N \mbox{
\,matrices \,}Q_{i,k} \mbox{ \,and \,}i\ge1.
$$
Using
$[\dt,(\dt)^i]=0$, we obtain that $Q_{i,k}=Q_i$ which
does not depend on $k$.
Note that since
$\dt=t^{-1}D$, we have $Q_1={\bf1}_N$ by (2.4).
\vs{5pt}

{\bf Lemma 2.3.} \
{\it
$P_1,P_2$ and
$Q_i-{\bf1}_N$ are strict
upper-triangular matrices for all $i\in2\Z_++1$.
}
\vs{5pt}

{\it Proof.}
So assume that $N>1$.
By (2.1) or (1.2),
we can deduce that
$$
\begin{array}{l}
-[i+1]_4(\dt)^{i-2}=
3[t^2\dt,[t^2\dt,(\dt)^i]]+2(2i-1)[t^3\dt,(\dt)^i],
\\[10pt]
0=[t^2\dt,[t^2\dt,[t^2\dt,(\dt)^i]]]
\\[4pt]
\phantom{0=}
+(i-1)(i-2)[t^4\dt,(\dt)^i]
+2(i-1)[t^2\dt,[t^3\dt,(\dt)^i]],
\\[10pt]
[i+1]_6(\dt)^{i-4}=
10[t^3\dt,[t^3\dt,(\dt)^i]-6(i-4)[t^5\dt,(\dt)^i]
\\[4pt]
\phantom{[i+1]_6(\dt)^{i-4}=}
-15[t^2\dt,[t^4\dt,(\dt)^i]],
\end{array}
$$
for $i\ge1$, where in general $[a]_j$ is a
notation similar to $[D]_j$ in (1.3)
(cf.~notation $[a]^j$ in (2.5)).
Here and below,
we make the convention that if a notion is not
defined but technically appears in an expression,
we always treat it as zero;
for instance,
$(\dt)^{i-2}=0$ if $i\le2$.
Applying these three
formulas to $Y_k$, we obtain
%%%%%%%%%%%%%%%%%%%%%%%
\\[10pt]
\hs{4pt}
$
-[i+1]_4Q_{i-2}=
3(P_{1,k-i+1}P_{1,k-i}Q_i
-2P_{1,k-i+1}Q_iP_{1,k}+Q_iP_{1,k+1}P_{1,k})
$\\[4pt]
\hs{4pt}
$
\phantom{-[i+1]_4Q_{i-2}=}
+2(2i-1)(P_{2,k-i}Q_i-Q_iP_{2,k}),
$\hfill(2.7)\\[10pt]
%%%%%%%%%%%%%%%%%%%%%%%
\hs{4pt}
$
0=P_{1,k-i+2}P_{1,k-i+1}P_{1,k-i}Q_i
-3P_{1,k-i+2}P_{1,k-i+1}Q_iP_{1,k}
$\\[4pt]
\hs{4pt}
$
\phantom{0=}
+3P_{1,k-i+2}Q_iP_{1,k+1}P_{1,k}
-Q_iP_{1,k+2}P_{1,k+1}P_{1,k}
$\\[4pt]
\hs{4pt}
$
\phantom{0=}
+(i-1)(i-2)(P_{3,k-i}Q_i-Q_iP_{3,k})
$\\[4pt]
\hs{4pt}
$
\phantom{0=}
+2(i-1)(P_{1,k-i+2}(P_{2,k-i}Q_i-Q_iP_{2,k})
$\\[4pt]
\hs{4pt}
$
\phantom{0=}
-(P_{2,k-i+1}Q_i-Q_iP_{2,k+1})P_{1,k}),
$\hfill(2.8)\\[10pt]
%%%%%%%%%%%%%%%%%%%%%%%
\hs{4pt}
$
[i+1]_6Q_{i-4}=10(P_{2,k-i+2}P_{2,k-i}Q_i
-2P_{2,k-i+2}Q_iP_{2,k} +Q_iP_{2,k+2}P_{2,k})
$\\[4pt]
\hs{4pt}
$
\phantom{[i+1]_6Q_{i-4}=}
-6(i-4)(P_{4,k-i}Q_i-Q_iP_{4,k})
$\\[4pt]
\hs{4pt}
$
\phantom{[i+1]_6Q_{i-4}=}
-15(P_{1,k-i+3}(P_{3,k-i}Q_i-Q_iP_{3,k})
$\\[4pt]
\hs{4pt}
$
\phantom{[i+1]_6Q_{i-4}=}
-(P_{3,k-i+1}Q_i-Q_iP_{3,k+1})P_{1,k}),
$\hfill(2.9)\\[10pt]
%%%%%%%%%%%%%%%%%%%%%%%%
for $i\ge1$.
We shall denote by $p_{i,k}^{(a,b)}$
the $(a,b)$-entry of the matrix $P_{i,k}$ and the
like for other matrices.
For a given position
$(a,b)$ with $1\le b\le a\le N$,
suppose inductively we have proved
$$
q_i^{(a_1,b_1)}=0
\eqno(2.10)
$$
for all $i\in2\Z_++1$
and for $a_1>a,\,b_1\le b$ or $a_1\ge a,\,b_1<b.$
Now for convenience, we denote
$$
p_{j,k}=p_{j,k}^{(a,a)},\;\;\;
p'_{j,k}=p_{j,k}^{(b,b)},\;\;\;
p_j=p_j^{(a,a)},\;\;\;p'_j=p_j^{(b,b)},\;\;\;
q_j=q_j^{(a,b)}
$$
for $j\in\Z.$
Assume that $i \in 2\Z_+ + 1$.
Using (2.10), by comparing
the $(a,b)$-entries in (2.7)-(2.9),
we obtain
%%%%%%%%%%%%%%%%%%%%%%%
\\[10pt]
\hs{4pt}
$
-[i+1]_4q_{i-2}=\bigl(3(p_{1,k-i+1}p_{1,k-i}
-2p_{1,k-i+1}p'_{1,k}+p'_{1,k+1}p'_{1,k})
$\\[4pt]
\hs{4pt}
$
\phantom{-[i+1]_4q_{i-2}=}
+2(2i-1)(p_{2,k-i}-p'_{2,k})\bigr)q_i,
$\hfill(2.11)\\[10pt]
%%%%%%%%%%%%%%%%%%%%%%%
\hs{4pt}
$
0=
 \bigl(p_{1,k-i+2}p_{1,k-i+1}p_{1,k-i}
-3p_{1,k-i+2}p_{1,k-i+1}p'_{1,k}
+3p_{1,k-i+2}p'_{1,k+1}p'_{1,k}
$\\[4pt]
\hs{4pt}
$
\phantom{0=}
-p'_{1,k+2}p'_{1,k+1}p'_{1,k}
+(i-1)(i-2)(p_{3,k-i}-p'_{3,k})
$\\[4pt]
\hs{4pt}
$
\phantom{0=}
+2(i-1)(
p_{1,k-i+2}(p_{2,k-i}-p'_{2,k})
$\\[4pt]
\hs{4pt}
$
\phantom{0=}
-(p_{2,k-i+1}-p'_{2,k+1})p'_{1,k})\bigr)q_i,
$\hfill(2.12)\\[10pt]
%%%%%%%%%%%%%%%%%%%%%%%
\hs{4pt}
$
[i+1]_6q_{i-4}=
\bigl(10(p_{2,k-i+2}p_{2,k-i}
-2p_{2,k-i+2}p'_{2,k}
+p'_{2,k+2}p'_{2,k})
$\\[4pt]
\hs{4pt}
$
\phantom{[i+1]_6q_{i-4}=}
-6(i-4)(p_{4,k-i}-p'_{4,k})
-15(p_{1,k-i+3}(p_{3,k-i}-p'_{3,k})
$\\[4pt]
\hs{4pt}
$
\phantom{[i+1]_6q_{i-4}=}
-(p_{3,k-i+1}-p'_{3,k+1})p'_{1,k})\bigr)q_i.
$\hfill(2.13)\\[10pt]
%%%%%%%%%%%%%%%%%%%%%%%
Applying
$[t^2\dt,t^{j+1}\dt]=(j-1)t^{i+2}\dt$ to
$y_k^{(b)}$ for $j=4,5$,
since $P_1,P_2$ are upper-triangular matrices,
using (2.5) and (2.6), we obtain
$$
\biggl\{
\begin{array}{l}
p_{4,k}=[\ol k]^5 + 10[\ol k]^3p_1
 + 10[\ol k]^2p_2  + 5\ol
kp_3 + p_4,
\\[6pt]
p_{5,k}=[\ol k]^6+15[\ol k]^4p_1
+20[\ol k]^3p_2+15[\ol k]^2p_3
+6\ol kp_4+p_5,
\end{array}
\eqno(2.14)
$$
where
$$
\begin{array}{l}
p_4 = -2(24p_1 + 12p_1^2 - 18p_2 + p_1p_2),\\[4pt]
p_5=5(-72p_1-34p_1^2+p_1^3+48p_2-6p_1p_2).
\end{array}
$$
We have similar formulas for $p'_{j,k},j=4,5$.
Applying $[t^3\dt,t^4\dt]=t^6\dt$
to $y_k^{(b)}$, we obtain the following
relation between $p_1$ and $p_2$, which is a
well-known relation for the Virasoro algebra
(cf.~[S1]).
$$
8p_1^2+4p_1^3-6p_1p_2+p_2^2 =0.
\eqno(2.15)
$$

First we make the following assumption
$$
q_i\ne0\mbox{ \ for some \ }i\in2\Z_++1.
\eqno(2.16)
$$
By
replacing $i$ by $i+2$ in (2.11), since
$[i+3]_4\ne0$ for $i\in2\Z_++1$, we see that
(2.16) holds for infinite many $i\in2\Z_++1$.
For
fixed $k$, we denote by $f_1(i),f_2(i),f_3(i)$
the coefficients of $q_i$ in (2.11)-(2.13)
respectively.
They are polynomials on $i$.
Then
(2.12) and (2.16) show that $f_2(i)=0$
for infinite many $i$.
Hence $f_2(i)=0$ for all
$i$.
Using (2.5) and (2.6) in (2.12), it is
straightforward to compute that the
coefficient of $i^4$ in $f_2(i)$ is $p_1-p'_1$.
Therefore,
$p'_1=p_1.$
Similarly, (2.11) and (2.13) show that
$$
g(i):=[i+1]_4[i-1]_4f_3(i)-[i+1]_6f_1(i-2)f_1(i)
$$
is zero for all $i$.
It is a little lengthy but
straightforward to compute that coefficient of
$i^{12}$ in $g(i)$ is $6p_1$ (using $p'_1=p_1$).
Thus $p_1=0$.
By (2.15), $p_2=0$.
Thus also
$p'_1=p'_2=0$.
Then (2.11) becomes
$[i+1]_4(q_i-q_{i-2})=0$.
From this we obtain that
$q_i=q_1$ for all $i\in2\Z_++1$.

Now we consider two cases:
First assume that $b<a$.
Then $q_1:=q_1^{(a,b)}=0$ (recall that
$Q_1={\bf1}_N$).
If (2.16) holds, then the above in
particular proves that $q_i=q_1=0$ for
all $i\in2\Z_++1$.
This contradicts (2.16).
Thus (2.16) cannot hold for any $i$, i.e.,
in this case we have $q_i=0$ for all
$i\in2\Z_++1$.

Next assume that $a=b$.
Then
$q_1:=q_1^{(a,a)}=1$ and so
(2.16) holds for at least $i=1$.
Thus the above
proves that $p_1=p_2=0,\,q_i=q_1$, i.e., in this
case we have $p_1^{(a,a)}=p_2^{(a,a)}=0$ and
$q_i^{(a,a)}=q_1^{(a,a)}=1$ for
all $i\in2\Z_++1$.

This proves the lemma.
\qed
\vs{5pt}

Lemma 2.3 shows that the diagonal elements
of $P_{j,k}$ are $[\ol k]^{j+1}$ for $j=1,2$, and
thus for all $j\ge1$ since $\Vir_+$ is generated
by $tD,t^2D$.
\vs{5pt}

{\bf Lemma 2.4.} \
{\it
$P_1=P_2=0$ and
$Q_i={\bf1}_N$ for all $i\in2\Z_++1$.
}
\vs{5pt}

{\it Proof.}
For a given position $(a,b)$
with $1\le a<b\le N$, suppose inductively
we have proved
$$
p_1^{(a_1,b_1)}=p_2^{(a_1,b_1)}=0,\;\;\;\;
q_i^{(a_1,b_1)}=\d_{a_1,b_1},
\eqno(2.17)
$$
for all $i\in2\Z_++1$ and for $a_1>a,b_1\le b$
or $a_1\ge a,b_1<b.$
Denote now
$$
p_{j,k}=p_{j,k}^{(a,a)},\;\;\; p_j=p_j^{(a,a)},
\;\;\;
p'_{j,k}=p_{j,k}^{(a,b)},\;\;\; p'_j=p_j^{(a,b)},
\;\;\;
q_i=q_i^{(a,b)},
$$
for
$j \in \Z,\,i \in 2\Z_++1$,
and denote
$$\wt
P_{j,k}= \left(\begin{array}{cc}p_{j,k}&p'_{j,k}\\
0&p_{j,k}\end{array}\right),
\ \
\wt P_j= \left(\begin{array}{cc}0&p'_j\\
0&0\end{array}\right),
\mbox{ \ \ and \ \ }
\wt Q_i= \left(\begin{array}{cc}1&q_i\\
0&1\end{array}\right).
$$
Then these
$2\times2$ matrices commute with each other.
By assumption (2.17), we see that (2.7)-(2.9)
still hold when we replace all matrices by
their corresponding matrices with tilde, and
we have similar formulas for $\wt
P_{j,k},j = 3,4,5$ as in (2.6) and (2.14)
(here now, $[\wt P_1,\wt P_2] = 0$).
Since
$\wt Q_i$ is invertible, from (2.8),
we obtain an equation on $\wt P_{i,k}$.
Using
(2.5) and (2.6) in this equation, we obtain that
$4[i]_3(3\wt P_1 - \wt
P_2) = 0$.
This shows that $\wt P_2 = 3\wt P_1$.
Then (2.7)
and (2.9) give
$$
\begin{array}{l}
[i+1]_4\wt Q_{i-2}=[i]_2(i^2-i+12\wt P_1-2)\wt Q_i,
\\[4pt]
[i+1]_6\wt Q_{i-4}=[i]_4(i^2-3i+30\wt P_1-4)\wt Q_i.
\end{array}
$$
Since $\wt Q_i$ are invertible, the above gives
$\wt P_1=0$ and so $\wt P_2=0$.
Then the above also
gives $\wt Q_i=\wt Q_1={\bf1}_2$ for $i\in2\Z_++1$.
This proves that (2.17) holds for $(a,b)$.
Thus we
have the lemma.
\qed
\vs{5pt}

Thus by Lemma 2.4 and (2.5), $P_{1,k}=[\ol k]^2$,
$P_{2,k}= [\ol k]^3$, $Q_i=1$, $i\in2\Z_++1$, are
all scalar matrices for $k>>0$.
By shifting the
grading index of $V_k$ if necessary, we can
suppose that $[\ol k]^2,[\ol k]^3\ne0$ and
(2.4) holds for $k\ge0$.
Applying
$[(\dt)^2,[(\dt)^2,t^2\dt]]=8(\dt)^3$ to $Y_0$,
we obtain that $Q_2^2=1$.
Thus by linear algebra,
$Q_2$ is a diagonalizable matrix.
Note that
$$
\si=(\dt)^2\cdot(t^3\dt)|_{V_0}
\eqno(2.18)
$$
is a linear transformation
on $V_0$ (recall (1.2) for the product ``$\cdot$''),
such that $\si Y_0=[\ol 2]^3Y_0Q_2$.
Thus by re-choosing the basis $Y_0$ and re-defining
$Y_k$ such that (2.4) holds for $k\ge0$ (then this
change of basis $Y_k$ does not effect
$P_{1,k},P_{2,k},Q_i,i\in2\Z_++1,$ since they are
scalar matrices), we can then suppose $Q_2$ is a
diagonal matrix (with the diagonal elements of
$Q_2$ being $\pm1$).
\vs{5pt}

{\bf Lemma 2.5.} \
{\it
For all $i,k\in\Z$
with $k,k+i\ge0$,
$P_{i,k}$ is a scalar matrix.
}
\vs{5pt}

{\it Proof.}
Using $[tD,t^{i-1}D] = (i - 2)t^iD$
and (2.5), by induction on $i$,
we obtain $P_{i,k} = [\ol k]^{i+1}$ for
$i \ge -1,\,k \ge 0$.
Thus assume that
$i = -i_1 \le -2,\,k + i \ge 0$.
Let $j$
be any integer such that $j > i_1$.
Applying
$(j + i_1)t^{j-i_1}D = [t^{-i_1}D,t^jD]$
to $Y_k$, we obtain
$$
(j + i_1)[\ol k]^{j-i_1+1}
=[\ol k]^{j+1}P_{-i_1,k+j}-[\ol
k - i_1]^{j+1}P_{-i_1,k}.
$$
By replacing $k$ by $k+j$
and replacing $j$ by $2j$, we obtain two other
equations respectively.
From these three equations,
one can easily deduce that $P_{-i_1,k}$ is a
scalar matrix.
\qed
\vs{5pt}

Since $W$ is generated by $\Vir\cup\{(\dt)^2\}$,
by induction on $j$, one can prove
$$
(t^{i+j}(\dt)^j)Y_k=Y_{k+i}P_{i,j,k}
\mbox{ \ for  some diagonal matrices \ }
P_{i,j,k},
\eqno(2.19)
$$
and for all
$i,j,k\in\Z$ with $j\ge1,\,k,i+k\ge0$.
\vs{5pt}

{\bf Lemma 2.6.} \
{\it
Denote by $V(a)$ the
$W$-submodule of $V$ generated by
$y_0^{(a)},\,a=1,...,N$.
Then $V(a)$ is a module
of the intermediate series such that
$V'=V(1)+...+V(N)$ is a direct sum of
$W$-submodules.
}
\vs{5pt}

{\it Proof.}
Since $U(W) = U(W_-)U(W_0 + W_+)$
and $V(a) = U(W)y_0^{(a)}$, by writing
$u \in  U(W)$ as a sum of
$u_1 \cdots  u_rw_1 \cdots  w_s$ for $u_i \in
 W_-$, $w_i \in  W_0 + W_+$, using (2.19), we
obtain by induction on $r + s$ that
${\rm dim\,}V(a)_k = 1$ for $k\ge0$.
Since $V(a)$ is also a $\Vir$-module, by [S2],
${\rm dim\,}V(a)_k = 1$ for all $k$ with
$k + \a \ne 0$.
Then by
(2.5) and the above lemmas, one can prove that
$V(a)$ is a subquotient module of $A_\a$ or
$B_\a$, i.e., $V(a)$ is a $W$-module of the
intermediate series (also cf.~[Z]).

For $a=1,...,N$, let $V'(a) =
 V(a) \cap \sum_{i\ne a}V(i)$.
Then obviously, $V'(a)_k = \{0\}$ for
$k \ge 0$.
Thus we must have $V'(a) = \{0\}$.
This proves the lemma.
\qed
\vs{5pt}

Now let $V''=V/V'$.
Then $V''$ is a
finite dimensional trivial module.
By induction
on the number $N+{\rm dim\,}V''$, one obtains
that $V$ is decomposable if $N\ge2$.
Thus $N=1$
and one can further deduce that $V$ is a module
of the intermediate series.
This proves Theorem
1.2(i).
\vs{5pt}

{\bf Corollary 2.7.} \
{\it
Suppose $V$ is
a uniformly bounded quasifinite $W$-module
satisfying $(2.2)$ and there exists $N\ge1$
such that ${\rm dim\,}V_i=N$ for all $i\in\Z$
with $\a+i\ne0$.
Fix $i_0\in\Z$ with $\a+i_0\ne0$
and fix a basis $Y_{i_0}$ of $V_{i_0}$.
Then there
exists a basis $Y_k$ of $V_k$ for all $k\in\Z$
with $\a+k\ne0$ such that
$(t^jD)Y_{i_0}=(\a+i_0)Y_{i_0+j}$ for all $j\in\Z$
with $\a+i_0+j\ne0$.
}
\qed
\vs{10pt}

\ni{\bf3. Quasifinite $\WW(\G,n)^{(1)}$-modules.}
Since Theorem 1.1(ii) is a special case of Theorem
1.2(ii), we shall prove Theorem 1.2(ii) (cf.~[S4]).
Thus assume that $\G$ is a group not isomorphic to
$\Z$ and $V$ is an indecomposable quasifinite weight
$\WW(\G,n)^{(1)}$-module such that there exists some
$\a=(\a_1,...,\a_n)\in\F^n$ with
(cf.~(2.2))
$$
V_\b=\{v\in V\,|\,D_iv=
(\a_i+\b_i)v,\,i=1,...,n\}\mbox{ \ for \ }\b\in\G.
$$
As the proof in [S4], $V$ is uniformly
bounded, and there exists $N\ge0$ such that
${\rm dim\,}V_\b=N$ for all $\b\in\G$ with
$\a+\b\ne0$.
For convenience, we shall now denote
$\ol\mu=\mu+\a$ for all $\mu \in \F^n$.

By [SZ2], we can suppose that all elements
$\g(i)=(\d_{1,i},...,\d_{n,i})$ for
$i=1,...,n,$ are in $\G$.
We
denote $\DD=\oplus_{i=1}^n\F D_i$ and define an
inner product on $\G\times\DD$ by
$$
\la\b,d\ra=\dsum{i=1}{n}\b_id_i\mbox{ \ \ for \ \ }
\b=(\b_1,...,\b_n)\in\G,\,d=\dsum{i=1}{n}d_iD_i
\in\DD.
\eqno(3.1)
$$
Then $\la\cdot,\cdot\ra$
is {\bf nondegenerate} in the sense that if
$\la\b,\DD\ra=0$ for some $\b\in\G$ then $\b=0$
and if $\la\G,d\ra=0$ for some $d\in\DD$
then $d=0$.

By (1.2) and (3.1), we have
$$
[t^\b d,t^\g
d']=t^{\b+\g}(\la\g,d\ra d'-\la\b, d'\ra d)
\mbox{ \ for \ }\b,\g\in\G,\;d,d'\in\DD.
\eqno(3.2)
$$
Fix an element $\g\in\G$
such that $\ol\g,\ol\g\pm\g(i),\ol\g\pm2\g(i)\ne0$
for $i=1,...,n$.
As in (2.18),
$$
\si_i=(t^{-2\g(i)}D_i(D_i - 1))\cdot(
t^{2\g(i)}D_i)|_{V_\g}\mbox{ \ for \ }i=1,...,n,
$$
are diagonalizable operators (note
that $(\dt)^2=t^2D(D - 1)$ and $t^3\dt=t^2D$
in (2.18)).
Since $\si_i,\,i=1,...,n$, commute with
each other, one can choose a basis $Y_\g$ of $V_\g$
such that $\si_i$ correspond to diagonal matrices.
Let $\b\in\G\bs\{0\}$ be any element such
that $\ol\g+\b\ne0$.
We shall define a basis
$Y_{\g+\b}$ of $V_{\g+\b}$ as follows:
One can
choose some $d\in\DD$ such that
$\la\ol\g,d\ra,\la\b,d\ra,\la\ol\g+\b,d\ra\ne0$.
Let $W(\b)={\rm span}\{t^{i\b}d^j\,|\,i\in\Z,j
\in\Z_+\bs\{0\}\}$ be a Lie
subalgebra of $\WW(\G,n)^{(1)}$, which is
isomorphic to $\WW(\Z,1)^{(1)}$ by (3.2)
(cf.~[S4]).
Denote
$V(\b)=\oplus_{i\in\Z}V_{\g+i\b}$.
Then $V(\b)$
is a uniformly bounded quasifinite $W(\b)$-module.
By Corollary 2.7, $t^\b
d|_{V_\g}:V_\g\rar V_{\g+\b}$ is bijective.
We
define $Y_{\g+\b}=
\la\ol\g+\b,d\ra^{-1}(t^\b d)Y_\g$.

Now as in
(2.19), one can prove by induction on
$|\mu| = \mu_1 + ... + \mu_n$ that $(t^\b
D^\mu)Y_\eta = Y_{\eta+\b}P_{\b,\mu,\eta}$ for
some diagonal matrices $P_{\b,\mu,\eta}$ and for
all $\b,\eta \in \G,\,\mu = (\mu_1,...,\mu_n)
 \in \Z_+^n\bs\{0\}$
with $\ol\eta,\ol\eta + \b \ne 0$.
Thus as
the proof of Lemma 2.6, we obtain that $V$ must
be a module of the intermediate series.
This
proves Theorem 1.2(ii).
\vs{10pt}

\cl{\bf References}
\vs{0pt}
\parskip .07 truein
\par\ni\hi3.5ex\ha1
{
\small
[BKLY] C.~Boyallian, V.~Kac, J.~Liberati and C.~Yan,
 Quasifinite highest weight modules of the Lie
 algebra of matrix differential operators on
 the circle, {\it J.~Math.~Phys.} {\bf39}
 (1998), 2910-2928.
 \par\ni\hi3.5ex\ha1
[C] V.~Chari, Integrable representations of
 affine Lie algebras, {\it Invent. Math.}
 {\bf85} (1986), 317-335.
 \par\ni\hi3.5ex\ha1
[FKRW] E.~Frenkel, V.~Kac, R.~Radul and W.~Wang,
 $\WW_{1+\infty}$ and $\WW(gl_N)$ with central
 charge $N$, {\it Comm.~Math.~Phys.} {\bf170}
 (1995), 337-357.
 \par\ni\hi3.5ex\ha1
[KL] V.~Kac and J.~Liberati, Unitary
 quasi-finite representations of $W_\infty$,
 {\it Lett.~Math.~Phys.} {\bf53} (2000), 11-27.
 \par\ni\hi3.5ex\ha1
[KP] V.~Kac and D.~Peterson, Spin and wedge
 representations of infinite dimensional Lie
 algebras and groups, {\it Proc.~Nat.~Acad.~Sci.~USA.}
 {\bf 78} (1981), 3308-3312.
 \par\ni\hi3.5ex\ha1
[KR1] V.~Kac and A.~Radul, Quasi-finite
 highest weight modules over the Lie algebra of
 differential operators on the circle,
 {\it Comm.~Math.~Phys.} {\bf157} (1993), 429-457.
 \par\ni\hi3.5ex\ha1
[KR2] V.~Kac and A.~Radul, Representation theory
 of the vertex algebra $\WW_{1+\infty}$,
 {\it Trans. Groups} {\bf1} (1996), 41-70.
 \par\ni\hi3.5ex\ha1
[KWY] V.~Kac, W.~Wang and C.~Yan, Quasifinite
 representations of classical Lie subalgebras
 of $\WW_{1+\infty}$, {\it Adv.~Math.}
 {\bf139} (1998), 46-140.
 \par\ni\hi3.5ex\ha1
[M] O.~Mathieu, Classification of Harish-Chandra
 modules over the Virasoro Lie algebra,
 {\it Invent. Math.} {\bf 107} (1992), 225-234.
 \par\ni\hi3.5ex\ha1
[S1] Y.~Su, A classification of indecomposable
 $sl_2(\C)$-modules and a conjecture of Kac on
 irreducible modules over the Virasoro algebra,
 {\it J.~Algebra} {\bf 161} (1993), 33-46.
 \par\ni\hi3.5ex\ha1
[S2] Y.~Su, Indecomposable modules over the
 Virasoro algebra, {\it Science in China A}
 {\bf44} (2001), 980-983.
 \par\ni\hi3.5ex\ha1
[S3] Y.~Su, 2-Cocycles on the Lie algebras of
 generalized differential operators,
 {\it Comm.~Algebra} {\bf30} (2002), 763-782.
 \par\ni\hi3.5ex\ha1
[S4] Y.~Su, Classification of quasifinite
 modules over the Lie algebras of Weyl type,
 {\it Adv.~Math.} {\bf 174} (2003), 57-68.
 \par\ni\hi3.5ex\ha1
[S5] Y.~Su, Classification of Harish-Chandra
 modules over the higher rank Virasoro algebras,
 {\it Comm. Math.~Phys.} {\bf240} (2003), 539-551.
 \par\ni\hi3.5ex\ha1
[S6] Y.~Su, Quasifinite representations of a
 Lie algebra of Block type, {\it J.~Algebra}
 {\bf276} (2004), 117-128.
 \par\ni\hi3.5ex\ha1
[SZ1] Y.~Su, K.~Zhao, Simple algebras of Weyl
 type, {\it Science in China A}
 {\bf44} (2001), 419-426.
 \par\ni\hi3.5ex\ha1
[SZ2] Y.~Su, K.~Zhao, Isomorphism classes and
 automorphism groups of algebras of Weyl type,
 {\it Science in China A} {\bf45} (2002), 953-963.
 \par\ni\hi3.5ex\ha1
[Z] K.~Zhao, The classification of a kind of
 irreducible Harish-Chandra modules over
 algebras of differential operators,
 {\it (Chinese) Acta Math. Sinica}
 {\bf37} (1994), 332-337.
}

\end{document}